\documentclass[11pt,twoside]{article}
\usepackage{amsmath,amsthm,amssymb,amscd,ascmac, amsfonts,fancybox}
\usepackage{mathrsfs}
\usepackage{url}
\allowdisplaybreaks[4]

\usepackage[dvips]{graphicx,color,psfrag}

\theoremstyle{plain}
\newtheorem{thm}{Theorem}

\newtheorem{lem}[thm]{Lemma}
\newtheorem{cor}[thm]{Corollary}

\newtheorem{prob}[thm]{Problem}

\theoremstyle{definition}

\newtheorem{rmk}[thm]{Remark}
\newtheorem{exa}[thm]{Example}

\begin{document}
\title{A generalization of a Ramanujan's exercise and Fibonacci polynomials}
\author{Genki Shibukawa}
\date{\empty}
\maketitle

\begin{abstract}
We give a generalization of a Ramanujan's exercise for high school students. 
Our results can be regarded as a variation of the factorization formula of $x^{n} - 1$. 
\end{abstract}

\section{Introduction}
Ramanujan pointed out the following elementary fact: {\it{the quadratic polynomial $at^{2}+t-a^{-1}$ is a factor of the polynomial $a^{5}t^{10}+11t^{5}-a^{-5}$}} (see \cite{TV} p9, Exercise 2). 
Thiruvenkatachar and Venkatachaliengar \cite{TV} generalized this fact as follows: 
\begin{thm}[Theorem~1.3.1 \cite{TV}]
\label{thm:Ramanujan's exercise}
The polynomial $at^{2}+t-a^{-1}$ is a factor of 
$$
a^{2n+1}t^{2(2n+1)}+Nt^{(2n+1)}-a^{-(2n+1)} \quad (n\geq 1),
$$
where
$$
N
   =
   (2n+1)+\sum_{r=1}^{n}\left\{\binom{n+r}{2r+1}+\binom{n+r+1}{2r+1}\right\}.
$$
\end{thm}
The Ramanujan's result is the particular case $n=2$ in Theorem~\ref{thm:Ramanujan's exercise}. 
Since the polynomial $at^{2}+t-a^{-1}=a^{-1}((at)^{2}+at-1)$ can be rewritten as $a^{-1}(x^{2}+x-1)$, Theorem~\ref{thm:Ramanujan's exercise} is equivalent to the following assertion:
\begin{align}
\label{eq:Ramanujan's exercise}
x^{2}+x-1 \mid x^{2(2n+1)}+Nx^{2n+1}-1.
\end{align}
We give a generalization of Theorem~\ref{thm:Ramanujan's exercise} or (\ref{eq:Ramanujan's exercise}).
\begin{thm}
\label{thm:main theorem}
For any non-negative integer $m\geq 0$, we define the polynomial $f_{m+1}(x;a,b) \in \mathbb{Z}[a,b][x]$ by 
$$
f_{m+1}(x;a,b)
   :=
   x^{2m+2}
   +\left[bF_{m}(a,b) - F_{m+2}(a,b)\right] x^{m+1}
   +b^{m+1},
$$
where $F_{m}(a,b)$ be the (bivariate) Fibonacci polynomial:
$$
F_{n+1}(a,b)
   :=
   \sum_{k=0}^{\lfloor \frac{n}{2}\rfloor}
   \binom{n-k}{k}a^{n-2k}(-b)^{k} \quad (n\geq 0), \quad 
   F_{0}(a,b)
   :=
   0.
$$
Then we have 
\begin{align}
\label{eq:main result}
f_{m+1}(x;a,b)
   &=
   (x^{2}-ax+b)
   \left[\sum_{j=0}^{m-1}F_{j+1}(a,b)x^{2m-j}
   +\sum_{j=0}^{m}b^{m-j}F_{j+1}(a,b)x^{j}\right] \\
\label{eq:main result2}
   &=
   x^{2m+2}-L_{m+1}(a,b)x^{m+1}+b^{m+1}.
\end{align}
Here the polynomial $L_{m}(a,b) \in \mathbb{Z}[a,b]$ is defined as the unique solution of the following difference equation:
$$
x_{m+2}-ax_{m+1}+bx_{m}=0, \quad x_{0}:=2, \quad x_{1}:=a.
$$
\end{thm}
By dividing (\ref{eq:main result}) by $x^{2}-ax+b$, we obtain
\begin{align}
\frac{f_{m+1}(x;a,b)}{x^{2}-ax+b}
   =
   \sum_{j=0}^{m-1}F_{j+1}(a,b)x^{2m-j}
   +\sum_{j=0}^{m}b^{m-j}F_{j+1}(a,b)x^{j}.
\end{align}
This formula is regarded as an analogue of the well-known geometric series formula:
$$
\frac{x^{m+1}-a^{m+1}}{x-a}
   =
   \sum_{j=0}^{m}x^{j}a^{m-j}.
$$

\section{Preliminaries}
In this section, we mention some fundamental properties of the Fibonacci polynomial $F_{n}(a,b)$ and others (see \cite{K} Chapter~43 for the detail). 
\begin{lem}
\label{lem:key lem}
The unique solution of the recursion:
$$
x_{n+2}-ax_{n+1}+bx_{n}=0, \quad x_{0}:=0, \quad x_{1}:=1
$$
is 
\begin{align}
x_{n}=F_{n}(a,b) \quad (n\geq 0).
\end{align}
In particular, we have
\begin{align}
\label{eq:Fibonacci}
F_{n}(\alpha +\beta ,\alpha \beta )
   =
   \frac{\alpha ^{n} - \beta ^{n}}{\alpha  - \beta }
   =
   \sum_{k = 0}^{n - 1}
      \alpha ^{n - 1 - k}\beta ^{k}.
\end{align}
\end{lem}
Under the following, we assume that $\alpha $ and $\beta $ are two roots of the quadratic polynomial:
$$
x^{2}-ax+b=(x-\alpha )(x-\beta ).
$$
Since $\alpha ^{n}$ and $\beta ^{n}$ satisfy the difference equation $x_{n+2}-ax_{n+1}+bx_{n}=0$ and 
$$
\alpha ^{0} + \beta ^{0} = 2, \quad \alpha ^{1} + \beta ^{1} = a,
$$
we obtain the following result.
\begin{lem}
\label{lem:key lem2}
For any non-negative integer $n$, we have
\begin{align}
\label{eq:Lucas}
L_{n}(\alpha + \beta ,\alpha \beta ) = \alpha ^{n} + \beta ^{n}.
\end{align}
\end{lem}

We define the polynomial $\Phi _{m+1}(x;a,b) \in \mathbb{Z}[a,b][x]$ by
\begin{align}
\Phi _{m+1}(x;a,b)
   :=
   \prod_{\substack{1\leq j\leq m+1 \\ \mathrm{gcd}(j,m+1)=1}}
   \left(x-\alpha e^{\frac{2\pi \sqrt{-1}j}{m+1}}\right)\left(x-\beta e^{\frac{2\pi \sqrt{-1}j}{m+1}}\right).
\end{align}
This polynomial $\Phi _{m+1}(x;a,b)$ is written as the product of two cyclotomic polynomials:
\begin{align}
\Phi _{m+1}(x;a,b)
   =
   b^{\varphi (m + 1)}
   \Phi _{m + 1}\left(\frac{x}{\alpha }\right)
   \Phi _{m + 1}\left(\frac{x}{\beta }\right),
\end{align}
where $\varphi (m + 1) = \mathrm{deg}(\Phi _{m + 1}(x))$ be the Euler totient function and $\Phi _{m + 1}(x)$ be the cyclotomic polynomial:
$$
\Phi _{m + 1}(x)
   :=
   \prod_{\substack{1\leq j\leq m+1 \\ \mathrm{gcd}(j,m+1)=1}}
   \left(x - e^{\frac{2\pi \sqrt{-1}j}{m+1}}\right).
$$

By the definition of $\Phi _{m+1}(x;a,b)$, the well-known factorization:
$$
x^{m+1}-a^{m+1}
   =
   \prod_{j=0}^{m}\left(x-e^{\frac{2\pi \sqrt{-1}j}{m+1}}a\right)
   =
   a^{m+1}\prod_{d \mid m+1}\Phi _{d}\left(\frac{x}{a}\right),
$$
and Lemma~\ref{lem:key lem2}, we obtain the following result. 
\begin{lem}
\label{lem:key lem3}
For any non-negative integer $m$, we have 
\begin{align}
\label{eq:cor:main result2-1}
x^{2m+2}-L_{m+1}(a,b)x^{m+1}+b^{m+1}
   &=
   \prod_{j=0}^{m}\left(x-\alpha e^{\frac{2\pi \sqrt{-1}j}{m+1}}\right)\left(x-\beta e^{\frac{2\pi \sqrt{-1}j}{m+1}}\right) \\
   &=
   \prod_{d \mid m+1}\Phi _{d}(x;a,b).
\end{align}
\end{lem}
These results are regarded as a quadratic analogue of the following well-known results:
$$
x^{m+1}-a^{m+1}
   =
   \prod_{j=0}^{m}\left(x-e^{\frac{2\pi \sqrt{-1}j}{m+1}}a\right)
   =
   a^{m+1}\prod_{d \mid m+1}\Phi _{d}\left(\frac{x}{a}\right).
$$

\section{Proof of Theorem~\ref{thm:main theorem}}
{\bf{Proof of Theorem~\ref{thm:main theorem}}} 
By direct calculus, we have (\ref{eq:main result}):
\begin{align*}
\text{LHS of (\ref{eq:main result})}
   &=
   x^{2m+2}+\sum_{j=1}^{m}(F_{j+1}(a,b)-aF_{j}(a,b)+bF_{j-1}(a,b))x^{2m+2-j} \nonumber \\
   & \quad +
   (2bF_{m}(a,b)-aF_{m+1}(a,b))x^{m+1} \\
   & \quad +
   b^{m+1}
   +\sum_{j=1}^{m-1}b^{m+1-j}(F_{j+1}(a,b)-aF_{j}(a,b)+bF_{j-1}(a,b))x^{j} \nonumber \\
   &=
   x^{2m+2} + (bF_{m}(a,b)-F_{m+2}(a,b))x^{m+1} + b^{m+1}.
\end{align*}
Here the second equality follows from Lemma~\ref{lem:key lem}.

From (\ref{eq:main result}), we have
$$
f_{m+1}(\alpha ;a,b)=f_{m+1}(\beta ;a,b)=0.
$$
By the definition of $f_{m+1}(x ;a,b)$, we have
$$
f_{m+1}\left(e^{\frac{2\pi \sqrt{-1}j}{m+1}}x;a,b \right)=f_{m+1}(x;a,b) \quad (j=0,1,\ldots, m).
$$
Then we give all the roots of the polynomial $f_{m+1}(x;a,b)$:
\begin{align}
\label{eq:cor:main result2-1}
f_{m+1}(x;a,b)
   =
   \prod_{j=0}^{m}\left(x-\alpha e^{\frac{2\pi \sqrt{-1}j}{m+1}}\right)\left(x-\beta e^{\frac{2\pi \sqrt{-1}j}{m+1}}\right).
\end{align}
Finally, by Lemma~\ref{lem:key lem3}, 
we obtain the conclusion (\ref{eq:main result2}). \qed 
\begin{rmk}
From the view point of symmetric functions, we obtain another proof of Theorem~\ref{thm:main theorem}. 
First, the formula (\ref{eq:main result2}) follows from (\ref{eq:Fibonacci}) and (\ref{eq:Lucas}):
\begin{align*}
\alpha \beta F_{m}(\alpha +\beta ,\alpha \beta )-F_{m+2}(\alpha +\beta ,\alpha \beta)
   &=
   \alpha \beta \frac{\alpha ^{m} - \beta ^{m}}{\alpha - \beta } - \frac{\alpha ^{m + 2} - \beta ^{m + 2}}{\alpha - \beta } \\
   &=
   - \alpha ^{m + 1} - \beta ^{m + 1} \\
   &=
   -L_{m + 1}(\alpha +\beta ,\alpha \beta ).
\end{align*}
We also obtain the formula (\ref{eq:main result}) straightforwardly by the following calculus:
\begin{align*}
& (x^{2}-(\alpha +\beta )x+\alpha \beta)
   \left[\sum_{j=0}^{m-1}F_{j+1}(\alpha +\beta,\alpha \beta)x^{2m-j}
   +\sum_{j=0}^{m}(\alpha \beta )^{m-j}F_{j+1}(\alpha +\beta,\alpha \beta)x^{j}\right] \\
& \quad =
   (x - \alpha )(x - \beta )
   \left[\sum_{j=0}^{m-1}\frac{\alpha ^{j + 1} - \beta ^{j + 1}}{\alpha -\beta }x^{2m-j}
   +\sum_{j=0}^{m}(\alpha \beta )^{m-j}\frac{\alpha ^{j + 1} - \beta ^{j + 1}}{\alpha -\beta }x^{j}\right] \\
& \quad =
   x^{2m}\frac{(x - \alpha )(x - \beta )}{\alpha - \beta }
   \left(
   \alpha \frac{1 - \alpha ^{m}/x^{m}}{1-\alpha /x}
   - \beta \frac{1 - \beta ^{m}/x^{m}}{1-\beta /x}
   \right) \\
   & \quad \quad +
   (\alpha \beta )^{m + 1}\frac{(x - \alpha )(x - \beta )}{\alpha - \beta }
   \left(
   \alpha \frac{1 - x^{m + 1}/\beta ^{m + 1}}{1 - x/\beta }
   - \beta \frac{1 - x^{m + 1}/\alpha  ^{m + 1}}{1 - x/\alpha }
   \right) \\
   &=
   \left(x^{2m + 2} - \frac{\alpha ^{m + 1} - \beta ^{m + 1}}{\alpha - \beta }x^{m + 2} + \alpha \beta \frac{\alpha ^{m} - \beta ^{m}}{\alpha  - \beta }x^{m + 1}\right) \\
   & \quad \quad +
   \left(\frac{\alpha ^{m + 1} - \beta ^{m + 1}}{\alpha - \beta }x^{m + 2} - \frac{\alpha ^{m + 2} - \beta ^{m + 2}}{\alpha  - \beta }x^{m + 1} + (\alpha \beta )^{m + 1}\right) \\
   &=
   f_{m+1}(x;\alpha +\beta ,\alpha \beta ).
\end{align*}
\end{rmk}

From Theorem~\ref{thm:main theorem} and Lemma~\ref{lem:key lem3}, we obtain the following more precisely factorization formula of $f_{m+1}(x;a,b)$. 
\begin{cor}
\label{cor:main result2}
We have
\begin{align}
f_{m+1}(x;a,b)
\label{eq:cor:main result2-2}
   &=
   \prod_{d \mid m+1}\Phi _{d}(x;a,b),
\end{align}
and 
\begin{align}
\sum_{j=0}^{m-1}F_{j+1}(a,b)x^{2m-j}
   +\sum_{j=0}^{m}b^{m-j}F_{j+1}(a,b)x^{j}
   &=
   \prod_{j=1}^{m}\left(x-\alpha e^{\frac{2\pi \sqrt{-1}j}{m+1}}\right)\left(x-\beta e^{\frac{2\pi \sqrt{-1}j}{m+1}}\right) \\
   &=
   \prod_{j=1}^{m}f_{1}\left(e^{\frac{2\pi \sqrt{-1}j}{m + 1}};a,b\right).
\end{align}
\end{cor}

\begin{exa}[$a=b=-1$]
In this case, $F_{n + 1}(-1,-1)$ and $L_{n + 1}(-1,-1)$ equal to the Fibonacci number 
$$
F_{n + 1}
   :=
   \sum_{r = 0}^{\lfloor \frac{n}{2} \rfloor}
      \binom{n - r}{r}, \quad F_{0}:=0,
$$
and Lucas number 
$$
L_{n + 1}
   :=
   F_{n} + F_{n + 2}
$$ 
essentially:
\begin{align*}
F_{n+1}(-1,-1)
   &=
   \sum_{k=0}^{\lfloor \frac{n}{2}\rfloor}
   \binom{n-k}{k}(-1)^{n-2k}
   =
   (-1)^{n}
   \sum_{k=0}^{\lfloor \frac{n}{2}\rfloor}
   \binom{n-k}{k}
   =
   (-1)^{n}F_{n + 1}, \\
L_{n+1}(-1,-1)
   &=
   (-1)^{n}(F_{n} + F_{n + 2})
   =
   (-1)^{n}L_{n + 1}.
\end{align*}
Thus, we have
\begin{align}
f_{m+1}(x;-1,-1)
   &=
   x^{2m + 2} + (-1)^{m}L_{m + 1}x^{m + 1} + (-1)^{m + 1} \\
   &=
   (x^{2} + x - 1)
   \left(
   \sum_{j=0}^{m - 1}
      (-1)^{j}F_{j+1}x^{2m - j}
   +(-1)^{m}\sum_{j=0}^{m}
      F_{j+1}x^{j}\right) \\
   &=
   \prod_{d \mid m+1}\Phi _{d}(x;-1,-1).
\end{align}

The case of $m = 2n$ is a generalization of Theorem~\ref{thm:Ramanujan's exercise}. 
In fact, by simple calculus, 
\begin{align*}
N
   &=
   \sum_{r=0}^{n - 1}
      \binom{2n - r - 1}{r}
   +\sum_{r=0}^{n}
   \binom{2n - r + 1}{r}
   =
   F_{2n} + F_{2n + 2}.
\end{align*}
Then we have
\begin{align}
x^{2(2n+1)}+Nx^{2n+1}-1
   &=
   f_{2n + 1}(x;-1,-1) \\
   &=
   (x^{2} + x - 1)
   \left(
   \sum_{j=0}^{2n - 1}
      (-1)^{j}F_{j+1}x^{4n-j}
   +\sum_{j=0}^{2n}
      F_{j+1}x^{j}\right) \\
   &=
   \prod_{d \mid 2n+1}\Phi _{d}(x;-1,-1).
\end{align}
The first few examples $\Phi _{n}(x;-1,-1)$ and $f_{n}(x;-1,-1)$ are given by 
\begin{align*}
\Phi _{1}(x;-1,-1)
   &=
   x^{2}+x-1, \\
\Phi _{2}(x;-1,-1)
   &=
   x^{2}-x-1, \\
\Phi _{3}(x;-1,-1)
   &=
   x^{4}-x^{3}+2x^{2}+x+1, \\
\Phi _{4}(x;-1,-1)
   &=   
   x^{4}+3x^{2}+1, \\
\Phi _{5}(x;-1,-1)
   &=
   x^{8} - x^{7} + 2x^{6} - 3x^{5} + 5x^{4} + 3x^{3} + 2x^{2} + x + 1 \\
   &=   
   (x^{4}+2x^{3}+4x^{2}+3x+1)(x^{4}-3x^{3}+4x^{2}-2x+1), \\
\Phi _{6}(x;-1,-1)
   &=
   x^{4}+x^{3}+2x^{2}-x+1,
\end{align*}
and 
\begin{align*}
f_{1}(x;-1,-1)
   &=
   x^{2}+x-1 \\
   &=
   \Phi _{1}(x;-1,-1), \\
f_{2}(x;-1,-1)
   &=
   x^{4}-3x^{2}+1 \\
   &=
   (x^{2} + x - 1)(x^{2} - x - 1) \\
   &=
   \Phi _{1}(x;-1,-1)\Phi _{2}(x;-1,-1), \\
f_{3}(x;-1,-1)
   &=
   x^{6}+4x^{3}-1 \\
   &=
   (x^{2} + x - 1)(x^{4} - x^{3} + 2x^{2} + x + 1) \\
   &=
   \Phi _{1}(x;-1,-1)\Phi _{3}(x;-1,-1), \\
f_{4}(x;-1,-1)
   &=
   x^{8}-7x^{4}+1 \\
   &=
   (x^{2} + x - 1)(x^{6} - x^{5} + 2x^{4} - 3x^{3} - 2x^{2} - x - 1) \\
   &=
   \Phi _{1}(x;-1,-1)\Phi _{2}(x;-1,-1)\Phi _{4}(x;-1,-1), \\
f_{5}(x;-1,-1)
   &=
   x^{10}+11x^{4}-1 \\
   &=
   (x^{2} + x - 1)(x^{8} - x^{7} + 2x^{6} - 3x^{5} + 5x^{4} + 3x^{3} + 2x^{2} + x + 1) \\
   &=
   \Phi _{1}(x;-1,-1)\Phi _{5}(x;-1,-1), \\
f_{6}(x;-1,-1)
   &=
   x^{12}-18x^{4}+1 \\
   &=
   (x^{2} + x - 1)(x^{10} - x^{9} + 2x^{8} - 3x^{7} + 5x^{6} - 8x^{5} - 5x^{4} - 3x^{3} - 2x^{2} - x - 1) \\
   &=
   \Phi _{1}(x;-1,-1)\Phi _{2}(x;-1,-1)\Phi _{3}(x;-1,-1)\Phi _{6}(x;-1,-1).
\end{align*}
\end{exa}

\section{A generalization}
Theorem~\ref{thm:main theorem} can be generalized as follows. 
We define the polynomial $f(x)$ by 
$$
f(x)
   :=
   \prod_{j=1}^{r}(x - \alpha _{j})
   =
   \sum_{k = 0}^{r}
      (-1)^{k}a_{k}x^{r - k}.
$$
\begin{prob}
Consider 
\begin{align*}
\frac{1}{f(x)}
   \prod_{j=1}^{r}(x^{m + 1} - \alpha _{j}^{m + 1})
   &=
   \prod_{k=1}^{m}f\left(e^{\frac{2\pi \sqrt{-1}k}{m + 1}}x\right)
   =
   \sum_{j = 0}^{r(m + 1)}
   c_{j}x^{r(m + 1) - j} \in \mathbb{Q}[a_{1},\ldots,a_{r}][x]
\end{align*}
and evaluate all coefficients $c_{j} \in \mathbb{Q}[a_{1},\ldots,a_{r}]$ explicitly. 
\end{prob}
The case of $r = 1$ is the geometric series and the case of $r = 2$ is Theorem~\ref{thm:main theorem}. 
Cases of $r \geq 3$ are difficult. 
Even in the following simple case, partial coefficients can be written by the higher order Fibonacci numbers \cite{S}:
$$
F_{m + 1}^{(3)}:\ldots, -25,-11,-5,-2,-1,0,0,1,1,3,4,9,\ldots
$$ 
and higher order Lucas numbers:
$$
L_{m + 1}^{(3)}:\ldots, 57,26,11,6,2,3,1,5,4,13,\ldots,
$$ 
but not all coefficients can be determined explicitly. 
\begin{exa}
We define $f(x)$ by 
$$
f(x) := x^{3} - x^{2} - 2x + 1 = \left[x + 2\cos{\left(\frac{2\pi }{7}\right)}\right]\left[x + 2\cos{\left(\frac{4\pi }{7}\right)}\right]\left[x + 2\cos{\left(\frac{6\pi }{7}\right)}\right]
$$
and put
\begin{align*}
f_{m + 1}(x)
   :=
   \left[x^{m + 1} - \left(-2\cos{\left(\frac{2\pi }{7}\right)}\right)^{m + 1}\right]
   \left[x^{m + 1} - \left(-2\cos{\left(\frac{4\pi }{7}\right)}\right)^{m + 1}\right]
   \left[x^{m + 1} - \left(-2\cos{\left(\frac{6\pi }{7}\right)}\right)^{m + 1}\right].
\end{align*}
The first few examples $f_{m + 1}(x)$ are
\begin{align*}
f_{1}(x)
   &=
   x^{3} - x^{2} - 2x + 1, \\
f_{2}(x)
   &=
   x^{6} - 5x^{4} + 6x^{2} - 1 \\
   &=
   f_{1}(x)(x^{3} + x^{2} - 2x - 1), \\
f_{3}(x)
   &=
   x^{9} - 4x^{6} - 11x^{3} + 1 \\
   &=
   f_{1}(x)(x^{6} + x^{5} + 3x^{4} + 5x^{2} + 2x + 1), \\
f_{4}(x)
   &=
   x^{12} - 13x^{8} + 26x^{4} - 1 \\
   &=
   f_{1}(x)(x^{9} + x^{8} + 3x^{7} + 4x^{6} - 4x^{5} + x^{4} - 11x^{3} - 5x^{2} - 2x - 1) \\
   &=
   f_{1}(x)(x^{3} + x^{2} - 2x - 1)(x^{6} + 5x^{4} + 6x^{2} + 1), \\
f_{5}(x)
   &=
   x^{15} - 16x^{10} - 57x^{5} + 1 \\
   &=
   f_{1}(x)(x^{12} + x^{11} + 3x^{10} + 4x^{9} + 9x^{8} - 2x^{7} + 12x^{6} - x^{5} + 25x^{4} + 11x^{3} + 5x^{2} + 2x + 1).
\end{align*}
\end{exa}

\section*{Acknowledgments}
This work was supported by JST CREST Grant Number JP19209317 and JSPS KAKENHI Grant Number 21K13808.


\bibliographystyle{amsplain}

\medskip
\begin{flushleft}
Genki Shibukawa \\
Department of Mathematics \\
Graduate School of Science \\
Kobe University \\
1-1, Rokkodai, Nada-ku, Kobe, 657-8501, JAPAN\\
E-mail: g-shibukawa@math.kobe-u.ac.jp
\end{flushleft}

\end{document}